\documentclass[reqno]{amsart}

\begin{document}

\newtheorem{theorem}{Theorem}[section]
\newtheorem{proposition}{Proposition}[section]
\newtheorem{lemma}{Lemma}[section]
\newtheorem{corollary}[theorem]{Corollary}
\newtheorem{example}[theorem]{Example}
\newtheorem{remark}{Remark}{}
\title[\hfil Principal solutions and irrationality in number theory ]
{Principal solutions of recurrence relations and irrationality questions in number theory } 
\author[Angelo B. Mingarelli \hfilneg]{Angelo B. Mingarelli}   
\address{School of Mathematics and Statistics\\ 
Carleton University, Ottawa, Ontario, Canada, K1S\, 5B6\\ and (present address) \\ Departamento de Matem\'aticas, Universidad de Las Palmas de Gran Canaria, 
Campus de Tafira Baja\\ 35017 Las Palmas de Gran Canaria, Spain.}
\email[A. B. Mingarelli]{amingare@math.carleton.ca, amingarelli@dma.ulpgc.es}

\date{February 14, 2007}
\thanks{This research is partially funded by a NSERC Canada Research Grant. Gratitude is also expressed to the Department of Mathematics of the University of Las Palmas in Gran Canaria for its hospitality during the author's stay there as a Visiting Professor during the two-year period, 2006-2008}
%\subjclass[2000]{34A30, 34C10}

\begin{abstract}
We apply the theory of disconjugate linear recurrence relations to the study of irrational quantities in number theory. In particular, for an irrational number associated with solutions of three-term linear recurrence relations we show that there exists a four-term linear recurrence relation whose solutions allow us to show that the number is a quadratic irrational if and only if the four-term recurrence relation has a principal solution of a certain type. The result is extended to higher order recurrence relations and a transcendence criterion can also be formulated in terms of these principal solutions. When applied to the situation of powers of $\zeta(3)$ it is not known whether the corresponding four term recurrence relation does or does not have such a principal solution, however the method does generate new series expansions of powers of $\zeta(3)$ and $\zeta(2)$ in terms of Ap\'{e}ry's now classic sequences. \end{abstract}

\maketitle

\section{Introduction} 

\noindent Of the methods used today to test for the irrationality of a given number we cite two separate approaches, one which seems to have overtaken the other recently. The first method is a direct consequence of Ap\'{e}ry's landmark paper \cite{rap}, which uses two independent solutions of a specific three-term recurrence relation (see \eqref{ap1} below) to generate a series of rationals whose limit at infinity is $\zeta(3)$.  Many new proofs and surveys of such arguments have appeared since, e.g., Beukers \cite{fb}, Nesterenko \cite{yu1}, Fischler \cite{sf}, Cohen \cite{hc2}, Murty \cite{rm}, Badea \cite{cb}, Zudilin \cite{wz}, to mention a few in a list that is far from exhaustive.\\

\noindent The idea and the methods used in  Ap\'{e}ry's work \cite{rap} were since developed and have produced results such as Andr\'{e}-Jeannin's proof of the irrationality of the inverse sum of the Fibonacci numbers \cite{aj}, along with a special inverse sum of Lucas numbers \cite{aj2}, and Zudilin's derivation \cite{wz1} of a three-term recurrence relation for which there exists two rational valued solutions whose quotients approach Catalan's constant. In addition we cite Zudilin's communication \cite{wz3} of a four-term recurrence relation (third order difference equation) for which there exists solutions whose quotients converge to $\zeta(5)$, but no irrationality results are derived.\\

\noindent Another approach involves considering the vector space $V$ over $\mathbb{Q}$ spanned by the numbers $1, \zeta(3), \zeta(5), \ldots, \zeta(2n+1)$. Using a criterion of Nesterenko \cite{yu2} on the linear independence of a finite number of reals, Rivoal \cite{tr} proved that $\dim V \geq c\log n$ for all sufficiently large $n$. It follows that the list $\zeta(3), \zeta(5), \ldots$ contains infinitely many irrationals. Rivoal complements this result in \cite{tr2} by showing that at least one of the numbers $\zeta(5), \zeta(7), \ldots, \zeta(21)$ is irrational. In the same vein, Zudilin \cite{wz2} shows that at least one of $\zeta(5), \zeta(7), \zeta(9), \zeta(11)$ is irrational. \\

\noindent In this work we  use the theory of disconjugate or non-oscillatory three-, four-, and n-term linear recurrence relations on the real line to problems in number theory; generally, to questions about the irrationality of various limits obtained via quotients of solutions at infinity and, in particular, to the irrationality and possible quadratic and higher algebraic irrationality of $\zeta(3)$ where $\zeta$ is the classic Riemann zeta function. We recall that this classic number is defined simply as
$$\zeta(3) = \sum_{n=1}^{\infty}\frac{1}{n^3}.$$ 

\noindent The underlying motivation here is two-fold. First, one can investigate the question of the irrationality of a given number $L$ say, by starting with an appropriate infinite series for $L$, associating to it a three-term recurrence relation (and so possibly a non-regular continued fraction expansion) whose form is dictated by the form of the series in question, finding an independent rational valued solution of said recurrence relation and, if conditions are right (cf. Theorem~\ref{th1} below), deduce the irrationality of $L$. We show that this abstract construction includes at the very least Ap\'{e}ry's classic proof \cite{rap} of the irrationality of  $\zeta(3)$. \\

\noindent Next, in our trying to determine whether or not $\zeta(3)$ is an algebraic irrational \cite{ghh}, we specifically address the question of whether $\zeta(3)$ is algebraic of degree two or more over $\mathbb{Q}$. Although we cannot answer this claim uniquevocally at this time, we present an equivalent criterion for the quadratic irrationality of $\zeta(3)$, or for that matter, any other irrational that can be approximated by the quotient of two solutions of an appropriate three-term recurrence relation. In the case of $\zeta(3)$ the equivalent criterion (Theorem~\ref{th3}) referred to is a function of the asymptotic behavior of solutions of a specific linear four-term disconjugate recurrence relation (Theorem~\ref{th2}, itself of independent interest) in which the products of the classic Ap\'{e}ry numbers play a prominent role, and whose general solution is actually known in advance. We obtain as a result, that appropriate products of the Ap\'{e}ry numbers satisfy a four-term recurrence relation, that is, \eqref{eqqq} below (indeed, given any $m\geq 2$ there exists an $(m+2)-$term recurrence relation for which these numbers play a basic role). However, the products of these Ap\'{e}ry numbers are not sufficient in themselves to give us the quadratic irrationality of $\zeta(3)$. Still, our results show that the quadratic irrationality of $\zeta(3)$ would imply the non-existence of linear combinations of appropriate products of Ap\'{e}ry sequences generating a principal solution of a certain type for this four-term linear recurrence relation. The converse is also true by our results but we cannot show that such linear combinations do not exist. Hence, we cannot answer at this time whether $\zeta(3)$ is a quadratic irrational.\\

\noindent We extend said criterion for quadratic irrationality of limits obtained by means of Ap\'{e}ry type constructions, or from continued fraction expansions to a criterion for algebraic irrationality (an irrational satisfying a polynomial equation of degree greater than two with rational coefficients) over $\mathbb{Q}$ (Theorem~\ref{th6}). It is then a simple matter to formulate a criterion for the transcendence of such limits. Loosely speaking, we show that an irrational number derived as the limit of a sequence of rationals associated with a basis for a linear three-term recurrence relation is transcendental if and only if there exists an infinite sequence of linear $m-$term recurrence relations, one for each $m\geq 2$, such that each one lacks a nontrivial rational valued solution with special asymptotics at infinity (cf., Theorem~\ref{degn}). Finally, motivated by the results on the four-term recurrences (Theorem~\ref{th2}), we present in the Appendix to this article accelerated series representations for $\zeta(3)^m$, for $m=2,3,4,5$, and similar series for $\zeta(2)^m$, where we display the cases $m=2,3$ only leaving the remaining cases as examples that can be formulated by the reader.\\

\section{Preliminary results} 

\noindent We present a series of lemmas useful in our later considerations.

\begin{lemma}\label{lem2} Let $A_n, c_n \in \mathbb{R}$, $n \in \mathbb{N}$, be two given infinite sequences such that the series
\begin{equation}
\label{eq4} 
\sum_{n=1}^{\infty}\frac{1}{c_{n-1}A_nA_{n-1}}
\end{equation}
converges absolutely.
 Then there exists a sequence $B_n$ satisfying \eqref{eq6} such that 
\begin{equation}
\label{eq5}
\lim_{n\to \infty} \frac{B_n}{A_n} = \frac{B_0}{A_0} + \alpha\,\sum_{n=1}^{\infty}\frac{1}{c_{n-1}A_nA_{n-1}},
\end{equation}
where $\alpha = c_0(A_0B_1-A_1B_0).$
\end{lemma}

\begin{proof}(Lemma~\ref{lem2}) For the given sequences $A_n, c_n$ define the sequence $b_n$ using \eqref{eq3} below:
\begin{equation}
\label{eq3} 
b_n = \{c_nA_{n+1}+c_{n-1}A_{n-1}\}/{A_n},  \quad n \geq 1.
\end{equation}
Then, by definition, the $A_n$ satisfy the three-term recurrence relation 
\begin{equation}
\label{eq1}
c_ny_{n+1}+c_{n-1}y_{n-1} - b_ny_n=0, \quad n\in \mathbb{N},
\end{equation}
with $y_0=A_0, y_1=A_1$. Choosing the values $B_0, B_1$ such that $\alpha \neq 0$, we solve the two-term recurrence relation 
\begin{equation}
\label{eq6}
A_{n-1}B_{n} - A_nB_{n-1} = \frac{\alpha}{c_{n-1}},  \quad n \geq 1.
\end{equation}
for a unique solution, $B_n$. Observe that these $B_n$ satisfy the same recurrence relation as the given $A_n.$ Since $A_nA_{n-1}\neq 0$ by hypothesis, dividing both sides of \eqref{eq6} by $A_nA_{n-1}$ gives \eqref{eq5} upon summation and passage to the limit as $n \to \infty$, since the resulting series on the left is a telescoping series.
\end{proof}

\begin{lemma}
\label{lem3}
Consider \eqref{eq1} where $c_n >  0$, {$b_n - c_n -c_{n-1} > 0$}, for every $n \geq n_0 \geq 1$, and $\sum_{n=1}^{\infty} 1/c_{n-1} < \infty.$ Let $A_m, B_m \in \mathbb{R}$, $m \geq 1$, be two linearly independent solutions of \eqref{eq1}.  If $0 \leq A_0 < A_1$, then 
\begin{equation}
\label{eq7}
L\equiv \lim_{m\to \infty} \frac{B_m}{A_m}
\end{equation}
exists and is finite.
\end{lemma}

\begin{proof} (Lemma~\ref{lem3})  Since $c_n > 0$, $b_n - c_n -c_{n-1} > 0$, for every $n \geq n_0$, the equation \eqref{eq1} is non-oscillatory at infinity \cite{ph1} or \cite{wtp}, that is, every solution $y_n$ has a constant sign for all sufficiently large $n$. From discrete Sturm theory we deduce that every solution of \eqref{eq1} has a finite number of nodes, \cite{wtp}. As a result, the solution $A_n$, may, if modified by a constant, be assumed to be positive for all sufficiently large $n$. Similarly, we may assume that $B_n > 0$ for all sufficiently large $n$. Thus, write $A_n > 0, B_n > 0$ for all $n  \geq N.$ Once again, from standard results in the theory of three-term recurrence relations, there holds the Wronskian identity \eqref{eq6} for these solutions. The proof of Lemma~\ref{lem2}, viz. \eqref{eq6},  yields the identity
\begin{equation}
\label{eq8}
\frac{B_n}{A_n} - \frac{B_{n-1}}{A_{n-1}}= \frac{\alpha}{c_{n-1}A_nA_{n-1}},
\end{equation}
for each $n \geq 1$. Summing both sides from $n=N+1$ to infinity, we deduce the existence of the limit $L$ in \eqref{eq7} (possibly infinite at this point) since the tail end of the series has only positive terms and the left side is telescoping. \\

We now show that the eventually positive solution $A_n$ is bounded away from zero for all sufficiently large $n$. This is basically a simple argument (see Olver and Sookne \cite{os} and Patula~[\cite{wtp}, Lemma 2] for early extensions). Indeed, the assumption $0 \leq A_0 <A_1$ actually implies that $A_n$ is increasing for all sufficiently large $n$. A simple induction argument provides the clue. Assuming that $A_{k-1} \leq A_k$ for all $1 \leq k \leq n$, $$A_{n+1} = \{b_nA_n - c_{n-1}A_{n-1}\}/c_n \geq A_n\{b_n - c_{n-1}\}/c_n > A_n,$$ since $b_n - c_n -c_{n-1} > 0$ for all large $n$. The result follows. \\

Now since $A_n$ is bounded away from zero for large $n$ and $\sum_{n=N+1}^{\infty} 1/c_{n-1} < \infty$ by hypothesis, it follows that the series 
$$ \sum_{n=N+1}^{\infty} \frac{1}{c_{n-1}A_nA_{n-1}} < \infty,$$
that is, $L$ in \eqref{eq7} is finite.
\end{proof}

\begin{remark} {\rm The limit of the sequence $A_n$ itself may be {\rm a priori} finite. For applications to irrationality proofs, we need that this sequence $A_n \to \infty$ with $n$. A sufficient condition for this is provided below.}
\end{remark} 

\begin{lemma}\label{lem3a}  {\rm (Olver and Sookne \cite{os}, Patula [\cite{wtp}, Lemma 2])} Let $c_n >  0$, 
\begin{equation}\label{bcc}
b_n - c_n -c_{n-1} > \varepsilon_n c_n,
\end{equation} 
for all sufficiently large $n$, where $\varepsilon_n > 0$, and $\sum_{n=1}^{\infty}\varepsilon_n$ diverges. Then every increasing solution $A_n$ of \eqref{eq1} grows without bound as $n \to \infty$.
\end{lemma}

\noindent  For the basic notion of disconjugacy in its simplest form we refer the reader to Patula~\cite{wtp} or Hartman~\cite{ph1}, for a more general formulation. For our purposes, \eqref{eq1} is a disconjugate recurrence relation on $[0,\infty)$ if every non-trivial solution $y_n$ has at most one sign change for all $n \in \mathbb{N}$. The following result is essentially a consequence of Lemma~\ref{lem3} and Lemma~\ref{lem3a}.

\begin{lemma}
\label{lemm}
 Let $c_n > 0$ in \eqref{eq1} satisfy $\sum_{n=1}^{\infty}1/c_{n-1} < \infty$. Let $b_n \in \mathbb{R}$ be such that $$b_n - c_n -c_{n-1}>0,$$ for $n \geq 1$. Then 
\begin{enumerate}
\item Equation \eqref{eq1} is a disconjugate three-term recurrence relation on $[0,\infty)$ 
\item There exists a solution $A_n$ with $A_n > 0$ for all $n\in \mathbb{N}$, $A_n$ increasing and such that for any other linearly independent solution $B_n$ we have the relation 
\begin{equation*}
\left | L - \frac{B_{m}}{A_{m}}\right | \leq \beta \frac{1}{A_{m}^2},
\end{equation*}
for some suitable constant $\beta$, for all sufficiently large $m$, where $L$ is the limit.
\item If, in addition, we have \eqref{bcc} satisfied for some sequence $\varepsilon_n > 0$ etc., then the solution $A_n$ in item {\rm (2)} grows without bound, that is, $A_n \to \infty$ as $n \to \infty$. 
\end{enumerate}
\end{lemma}

\noindent Item (2) of the preceding lemma is recognizable by anyone working with continued fractions, \cite{ghh}. Of course, continued fractions have convergents (such as $A_n, B_n$ above) that satisfy linear three-term recurrence relations and their quotients, when they converge, converge to the particular number (here represented by $L$)  represented by the continued fraction.  In this article we view the limits of these quotients in terms of asymptotics of solutions of disconjugate  recurrence relations, with a particular emphasis on principal solutions.

\section{Main results}%%%%%%%%%%%%%%%%%%%%%%%%%%%%%

\begin{theorem} 
\label{th1}
Consider the three-term recurrence relation \eqref{eq1} where $b_n \in \mathbb{R}$, $c_n >  0$ for every $n \geq n_0 \geq 1$, and   the leading term $c_n$ satisfies
\begin{equation}
\label{cnn1}
\sum_{n=1}^{\infty} \frac{1}{c_{n-1}^{k}} <\infty,
\end{equation} 
for some $k \geq 1$. In addition, let \eqref{bcc} be satisfied for some sequence $\varepsilon_n > 0$, with $\sum_{n=1}^{\infty}\varepsilon_n = \infty$. 

Let $0 \leq A_0 < A_1$ be given and the resulting solution $A_m$ of \eqref{eq1}, satisfy $A_m \in \mathbb{Q}^+$ for all large $m$, and
\begin{equation}
\label{eq9a}
\sum_{n=1}^{\infty}\frac{1}{A_{n}^{\delta}} < \infty,
\end{equation}
for some $\delta$, where $0 < \delta< k^{\prime}$ and $k^{\prime}=k/(k-1)$ whenever $k>1$.

Next, let $B_m$ be a linearly independent solution such that $B_m \in \mathbb{Q}$ for all sufficiently large $m$ and such that for some sequences $d_m, e_m \in \mathbb{Z}^+$, we have $d_mA_m \in \mathbb{Z}^+$ and $e_mB_m \in \mathbb{Z}^+$, for all sufficiently large $m$, and
\begin{equation}
\label{eq9b}
\lim_{m \to \infty} \frac{{\rm lcm}\{d_m, e_m\}}{{A_m}^{1-\delta/k^{\prime}}} =0.
\end{equation}
Then $L$, defined in \eqref{eq7}, is an irrational number.
\end{theorem} 

\begin{proof} (Theorem~\ref{th1}) We separate the cases $k=1$ from $k>1$ as is usual in this kind of argument. Let $k=1$. With $A_n, B_n$ as defined, a simple application of Lemma~\ref{lem3} (see \eqref{eq8}) gives us that, for $m \geq N$,
\begin{equation}
\label{eq10}
\sum_{n=m+1}^{\infty} \left \{\frac{B_n}{A_n} - \frac{B_{n-1}}{A_{n-1}}\right \} = \alpha \sum_{n=m+1}^{\infty}\frac{1}{c_{n-1}A_nA_{n-1}},
\end{equation}
i.e., 
\begin{equation}
\label{eq11}
L - \frac{B_{m}}{A_{m}} = \alpha \sum_{n=m+1}^{\infty}\frac{1}{c_{n-1}A_nA_{n-1}}.
\end{equation}
Since $A_n$ is increasing for all $n \geq N$ (by Lemma~\ref{lem3}) we have $A_nA_{n-1} > A_{n-1}^2$ for such $n$. In fact, we also have $A_n \to \infty$ (by Lemma~\ref{lem3a}). Estimating \eqref{eq11} in this way we get
\begin{equation}\label{lbm}
\left | L - \frac{B_{m}}{A_{m}}\right | \leq  \alpha \sum_{n=m+1}^{\infty}\frac{1}{c_{n-1} A_{n-1}^2},
\end{equation} 
and since $A_{k}^2 > A_{m}^2$ for $k> m$ we obtain
\begin{equation}
\label{eq12}
\left | L - \frac{B_{m}}{A_{m}}\right | \leq \beta \frac{1}{A_{m}^2},
\end{equation}
where $\beta = \alpha \sum_{n=m+1}^{\infty}1/c_{n-1} < \infty$ is a constant. The remaining argument is conventional. 
Multiplying \eqref{eq12} by ${\rm lcm}\{d_m, e_m\}\cdot A_m$ for all large $m$, we find 
\begin{equation}
\label{eq13}
\left |{\rm lcm}\{d_m, e_m\} A_m L - B_m{\rm lcm}\{d_m, e_m\} \right | \leq \beta \frac{{\rm lcm}\{d_m, e_m\}}{A_{m}}.
\end{equation}
Assuming that $L = C/D$ is rational, where $C, D$ are relatively prime, we get 
\begin{equation*}
\left |{\rm lcm}\{d_m, e_m\} A_m C - B_m D {\rm lcm}\{d_m, e_m\} \right | \leq \beta D \frac{{\rm lcm}\{d_m, e_m\}}{A_{m}}.
\end{equation*}
But the left hand side is a non-zero integer for every $m$ (see \eqref{eq11}), while the right side goes to zero as $m \to \infty$ by \eqref{eq9b} with $k^{\prime}=\infty$. Hence it must eventually be less than $1$, for all large $m$, which leads to a contradiction. This completes the proof in the case $k=1$. \\

Let $k > 1$. We proceed as in the case $k=1$ up to \eqref{lbm}. That the solution $A_n$ as defined is increasing is a consequence of the proof of Lemma~\ref{lem3}. The fact that this $A_n \to \infty$ as $n \to \infty$ follows from Lemma~\ref{lem3a}. The existence of the limit $L$ is clear since the series consists of positive terms for all sufficiently large $m$. In order to prove that this $L$ is indeed finite we observe that 
$$\left | L - \frac{B_{m}}{A_{m}}\right | \leq \beta \left \{ \sum_{n=m+1}^{\infty}\frac{1}{A_n^{k^{\prime}}A_{n-1}^{k^{\prime}}}\right \}^{1/k^{\prime}}$$ where $\beta = \alpha \{\sum_{n=m+1}^{\infty}{1/c_{n-1}^{k}}\}^{1/{k}} <\infty$, by \eqref{cnn1}. Next, $A_n^{k^{\prime}}A_{n-1}^{k^{\prime}}$ $= A_n^\delta A_n^{k^{\prime}-\delta}A_{n-1}^{k^{\prime}}$ $\geq A_n^\delta A_m^{2k^{\prime}-\delta}$, for all sufficiently large $n$. Hence
\begin{equation}
\label{eq12a}
\left | L - \frac{B_{m}}{A_{m}}\right | \leq \beta^{\prime} \frac{1}{A_m^{2-\delta/k^{\prime}}},
\end{equation}
where $\beta^{\prime}  = \beta \{ \sum_{n=m+1}^{\infty}1/A_{n}^{\delta}\}^{1/k^{\prime}} < \infty$ by \eqref{eq9a}. Since $0 < \delta < k^{\prime}$, we get that $L$ is finite. Equation \eqref{eq12a} corresponds to \eqref{eq12} above. Continuing as in the case $k=1$ with minor changes, we see that \eqref{eq9b} is sufficient for the irrationality of $L$.\\ 
\end{proof}

\begin{remark}{\rm  Condition \eqref{eq9a} is not needed in the case $k=1$. This same condition is verified for corresponding solutions of recurrence relations of the form \eqref{eq1} with $c_{n}=n+1$, $b_n=an+b$ where $a>2$, for all sufficiently large indices. Note that the case $a=2$ is a borderline case. For example, for $a=2, b=1$, there are both bounded nonoscillatory solutions (e.g., $y_n=1$) and unbounded nonoscillatory solutions (e.g., $y_n = 1 + 3\Psi (n+1) +3\gamma$, where $\Psi(x) = (\log \Gamma(x))^{\prime}$ is the {\rm digamma function} and $\gamma$ is Euler's constant). Thus, for every pair of such solutions, the limit $L$ is either infinite or a rational number. For $a < 2$ all solutions are oscillatory, that is $y_ny_{n-1} < 0$ for arbitrarily large indices. Such oscillatory cases could also be of interest for number theoretical questions, especially so if the ratio of two independent solutions is of one sign for all sufficiently large $n$ (as in {\rm Zudilin}~\cite{wz3}).}
\end{remark}

\subsection{Consequences and discussions} 

The simplest consequences involve yet another interpretation of the proof of the irrationality of $\zeta(3)$ (and of $\zeta(2)$). It mimics many of the known proofs yet a large part of it involves only the theory of disconjugate three-term recurrence relations. Since the proofs are similar we sketch the proof for the case of $\zeta(3)$.\\

\begin{proposition} $\zeta (3)$ is irrational.
\label{prop1}
\end{proposition}

\begin{proof} (Proposition~\ref{prop1}) (originally due to Ap\'{e}ry \cite{rap}, cf., also Van der Poorten \cite{vdp}, Beukers \cite{fb}, Cohen \cite{hc1}). \\

\noindent Consider \eqref{eq1} with $c_n=(n+1)^3$ and $b_n = 34n^3+51n^2+27n+5$, $n\geq 0$. This gives the recurrence relation of Ap\'{e}ry, 
\begin{equation}\label{ap1}
(n+1)^3y_{n+1}+n^3y_{n-1} = (34n^3+51n^2+27n+5)y_n, \quad n\geq 1.
\end{equation}
\noindent Define two independent solutions $A_n, B_n$ of \eqref{ap1} by the initial conditions $A_0 =1, A_1=5$ and $B_0=0, B_{1}=6$. Then $b_n-c_n-c_{n-1} > 0$ for every $n\geq 0$. Since $0 < A_0 < A_1$ the sequence $A_n$ is increasing by Lemma~\ref{lem3} and tends to infinity with $n$. Note that, in addition to $c_n > 0$, \eqref{bcc} is satisfied for every $n\geq 1$ and $\varepsilon_n=1/n$, say. Hence, \eqref{ap1} is a disconjugate three-term recurrence relation on $[0,\infty)$. An application of Lemma~\ref{lem3} and Lemma~\ref{lem3a} shows that 
\begin{equation}
\label{eq14}
L\equiv \lim_{m\to \infty} \frac{B_m}{A_m}
\end{equation}
\noindent exists and is finite and, as a by-product, we get \eqref{eq5}, that is (since $B_0=0$), 
\begin{equation}
\label{lal}
L =  \alpha\,\sum_{n=1}^{\infty}\frac{1}{c_{n-1}A_nA_{n-1}},
\end{equation}
\noindent where $\alpha=6$ in this case. \\

\noindent Define non-negative sequences $A_n, B_n$ by setting 
\begin{equation}
\label{an}A_n = \sum_{k=0}^{n}{\binom{n}{k}}^2{\binom{n+k}{k}}^2,
\end{equation} 
and 
\begin{equation}
\label{bn}B_n = \sum_{k=0}^{n}{\binom{n}{k}}^2{\binom{n+k}{k}}^2\left\{\sum_{m=1}^{n} \frac{1}{m^3} + \sum_{m=1}^{k} \frac{(-1)^{m-1}}{2m^3\binom{n}{m}\binom{n+m}{m}} \right \}.
\end{equation} 
A long and tedious calculation (see Cohen~\cite{hc1}) gives that these sequences satisfy \eqref{ap1}, and thus must agree with our solutions (bearing the same name) since their initial values agree. That $L=\zeta (3)$ in \eqref{eq14} is shown directly by using these expressions for $A_n, B_n$. In addition, it is clear that $A_n \in \mathbb{Z}^+$ (so $d_n=1$ in Theorem~\ref{th1}) while the $B_n\in \mathbb{Q}^+$ have the property that if $e_m = 2{\rm lcm}[1,2,...n]^3$ then $e_mB_m \in \mathbb{Z}^+$, for every $m\geq 1$ (cf., e.g., \cite{vdp}, \cite{hc1} among many other such proofs). Hence the remaining conditions of Theorem~\ref{th1} are satisfied, for $k=1$ there. So, since it is known that asymptotically $e_m/A_m \to 0$ as $m \to \infty$ (e.g., \cite{vdp}), the result follows from said Theorem.
\end{proof}

\begin{remark}{\rm Strictly speaking, the number thoeretical part only comes into play after \eqref{lal}. If we knew somehow that the series in \eqref{lal} summed to $\zeta(3)$ independently of the relations \eqref{an}, \eqref{bn} that follow, we would have a more natural proof. This is not a simpler proof of the irrationality of $\zeta (3)$; it is simply a restatement of the result in terms of the general theory of recurrence relations, in yet another approach to the problem of irrationality proofs. The proof presented is basically a modification of Cohen's argument in \cite{hc1} recast as a result in the asymptotic theory of three-term recurrence relations. We also observe that a consequence of the proof is that {\rm ([Fischler~\cite{sf}, Remarque 1.3, p. 910-04])}, 
\begin{equation}
\label{eq15}
\zeta (3) =  6\,\sum_{n=1}^{\infty}\frac{1}{n^3A_nA_{n-1}},
\end{equation}
an infinite series that converges much faster (series acceleration) to $\zeta(3)$ than the original series considered by Ap\'{e}ry, that is 
\begin{equation}\label{eq16}
\zeta(3) = \frac{5}{2}\sum_{n=1}^{\infty} \frac{(-1)^{n-1}}{n^3\binom{2n}{n}}.
\end{equation}
For example, the first 5 terms of the series \eqref{eq15} gives 18 correct decimal places to $\zeta(3)$ while \eqref{eq16} only gives 4. At the end of this paper we provide some series acceleration for arbitrary integral powers of $\zeta(3)$.} \end{remark}

\noindent The preceding remark leads to the following natural scenario. Let's say that we start with the infinite series  
\begin{equation}
\label{caa} L =  \sum_{n=1}^{\infty} \frac{1}{n^3A_nA_{n-1}}
\end{equation}
where the terms $A_n$ are the Ap\'{e}ry numbers defined in \eqref{an} and the series \eqref{caa} has been shown to be convergent using direct means that is, avoiding the use of the recurrence relation \eqref{ap1}. Then, by Lemma~\ref{lem2} there exists a rational valued sequence $B_n$ such that both $A_n, B_n$ are linearly independent solutions of a three-term recurrence relation of the form \eqref{eq1}. The new sequence $B_n$ thus obtained must be a constant multiple of their original counterpart in \eqref{bn}. Solving for the $b_n$ using \eqref{eq5} would necessarily give the cubic polynomial in \eqref{ap1}, which has since been a mystery. Once we have the actual recurrence relation in question we can then attempt an irrationality proof of the number $L$ using the methods described, the only impediment being how to show that $e_mB_m \in \mathbb{Z}^+$ without having an explicit formula like \eqref{bn}.\\

\noindent The method can be summed up generally as follows: We start with an infinite series of the form
\begin{equation}
\label{caa1} L = \sum_{n=1}^{\infty} \frac{1}{c_{n-1}A_nA_{n-1}}
\end{equation}
where the terms $c_n , A_n \in \mathbb{Z}^+$, and the series \eqref{caa1} has been shown to be convergent to $L$ using some direct means. Then, by Lemma~\ref{lem2} there exists a rational valued sequence $B_n$ such that both $A_n, B_n$ are linearly independent solutions of \eqref{eq1} where the $b_n$, defined by \eqref{eq5} are rational for every $n$. If, in addition, we have for example,
\begin{equation}\label{1cn}
\sum_{n=1}^{\infty} 1/c_{n-1} < \infty,
\end{equation} 
along with \eqref{bcc} we can then hope to be in a position so as to apply Theorem~\ref{th1} and obtain the irrationality of the real number $L$. Of course, this all depends on the interplay between the growth of the $d_nA_n$ at infinity and the rate of growth of the sequence $e_nB_n$ required by said Theorem (see \eqref{eq9b}). The point is that the relation \eqref{eq12} used by some to obtain irrationality proofs for the number $L$, is actually a consequence of the theory of disconjugate three-term recurrence relations. In fact, underlying all this is Lemma~\ref{lemm}. \\

\noindent The next two results are expected and included because their proofs are instructive for later use.
\begin{proposition} \label{prop2} The only solution of \eqref{ap1} whose values are all positive integers is, up to a constant multiple, the solution $A_n$ in \eqref{an}.
\end{proposition}

\begin{proof} (Proposition~\ref{prop2}) If possible, let $D_n$ be another integer valued solution of \eqref{ap1}. Then $D_n=aA_n+bU_n$, for every $n\in \mathbb{N}$ where $a,b\in \mathbb{R}$ are constants. Using the initial values $A_0=1, A_1=5$, $U_0=\zeta(3), U_1=5\zeta(3)-6$, in the definition of $D_n$, we deduce that $a=D_0-(5D_0-D_1)\zeta(3)/6$ and $b=(5D_0-D_1)/6.$ Thus,
$$D_n = D_0 A_n - (5D_0-D_1)B_n/6, \quad n \geq 1,$$ where the coefficients of $A_n, B_n$ above are rational numbers. By hypothesis, the sequence $D_n$, $n\in \mathbb{N}$ is integer valued. But so is $A_n$; thus $D_n-D_0A_n \in \mathbb{Z}$ for all $n$. Therefore, for $5D_0-D_1\neq 0$, we must have that $(5D_0-D_1)B_n/6 \in \mathbb{Z}$ for all $n$, which is impossible for sufficiently large $n$ (see \eqref{bn}). Hence $5D_0-D_1= 0$, and this shows that $D_n$ must be a multiple of $A_n$.
\end{proof}

\begin{proposition}
\label{prop3} The solution $B_n$ of Ap\'{ery} is not unique. That is there exists an independent strictly rational (i.e, non-integral) solution $D_n$ of \eqref{ap1} such that $$\frac{1}{3}{\rm lcm}[1,2,...,n]^3 D_n \in \mathbb{Z}^+$$ for all $n$.
\end{proposition}
\begin{proof} (Proposition~\ref{prop3}) A careful examination of the proof of Proposition~\ref{prop2} shows that the solution $B_n$ defined in \eqref{bn} is not the only solution of \eqref{ap1} with the property that $2{\rm lcm}[1,2,...,n]^3 B_n \in \mathbb{Z}^+$ for all $n$. Indeed, the solution $D_n$, defined by setting $D_0 = 1$, $D_1=1$ and $D_n = D_0 A_n - (5D_0-D_1)B_n/6$, for $n \geq 1,$ has the additional property that $2{\rm lcm}[1,2,...,n]^3 D_n/6 \in \mathbb{Z}^+$ for all $n$. Thus, the claim is that the quantity $2{\rm lcm}[1,2,...,n]^3 B_n$ is always additionally divisible by 6, for every $n \in \mathbb{N}$. That is, it suffices to show that ${\rm lcm}[1,2,...,n]^3 B_n$ is divisible by 3. But this can be accomplished by considering the contribution of this additional divisor to the p-adic valuation, $v_p(\cdot)$, of one term of the third sum in \eqref{bn}. Consider Cohen's proof [\cite{hc1}, Proposition 3] that $2{\rm lcm}[1,2,...,n]^3 B_n \in \mathbb{Z}^+$ for all $n$. There he shows that the quantity 
$$v= v_p\left (\frac{d_n^3\binom{n+k}{k}}{m^3\binom{n}{m}\binom{n+m}{m}}\right ) \geq \ldots \geq (v_p(d_n)-v_p(m))+(v_p(d_n)-v_p(d_k)) \geq 0,$$
where $d_n = {\rm lcm}[1,2,...,n]^3.$ Observe that insertion of the factor $1/3$ on the left only decreases the right side by 1 for the 3-adic valuation (see [\cite{hc1}, p.VI.5],) and then, keeping track of the other two terms above on the right and the fact that they are not zero we see that the inequality is still valid. Of course, one cannot do better than the divisor `6' in this respect since $B_1=6$. 
\end{proof}

\begin{remark}{\rm  A simple heuristic argument in the case of $\zeta(5)$ shows that if we are looking for recurrence relations of the form \eqref{eq1} with $c_n=(n+1)^5$ and $b_n$ some quintic polynomial in $n$, and we want an integral-valued solution other than the trivial ones, then we must have the coefficient of the leading term of the quintic superior to $150$ in order for the asymptotics to work out at all. The subsequent existence of a second solution $S_n$ with the property that $c\cdot{\rm lcm}[1,2,...,n]^5 S_n \in \mathbb{Z}^+$ for all $n$, where $c$ is a universal constant, is then not out of the question and could lead to an irrationality proof of this number. However, it is not at all clear to us that such a (non-trivial) quintic exists.}
\end{remark}

\noindent The basic advantage of the formalism of recurrence relations lies in that every element in $\mathbb{Q}(\zeta(3))$ can be approximated by `good' rationals, that is appropriate linear combinations of the $A_n, B_n$ in \eqref{an}, \eqref{bn}. For example, the series considered by Wilf \cite{hsw}
$$\sum_{n=1}^{\infty} \frac{1}{n^3(n+1)^3(n+2)^3} = \frac{29}{32} - \frac{3}{4}\zeta(3),$$
derived as a result of the use of the WZ algorithm e.g., \cite{tad}, has a counterpart via \eqref{ap1}. The solution $C_n$ of \eqref{ap1} defined by $C_n=(29/32)A_n-(3/4)B_n$ has the property that $C_n/A_n \to (29/32)-(3/4)\zeta(3)$ as $n \to \infty$, and the convergence of these fractions is sufficiently rapid as to ensure the irrationality of its limit, but this does not appear to be so for Wilf's series, even though it is an `accelerated' series. A similar comment applies to the series
$$\sum_{n=1}^{\infty} \frac{1}{(n+1)^3(n+2)^3(n+3)^3(n+4)^3(n+5)^3} = \frac{5}{768}\zeta(3) -\frac{10385}{98304},$$ also derived in \cite{hsw}. We point out that the above two series can also be summed more simply by using the method of partial fractions.\\

\noindent  The usefulness of so-called {\it dominant} and {\it recessive} solutions in the theory (also called {\it principal} solutions by some) is apparent in the following discussion regarding the overall nature of the solutions of \eqref{ap1}. As noted earlier, $A_n > 0$ for every $n$, $A_n$ is increasing, and the series in \eqref{eq5} converges.  In addition, by defining the solution $U_n= \zeta(3)A_n - B_n$, we see that $U_n/A_n\to 0$ as $n \to \infty$ (see the proof of Proposition~\ref{prop1}). Hence, by definition, $A_n$ (resp. $U_n$) is a dominant (resp. recessive) solution of the disconjugate equation \eqref{ap1} on $[0,\infty)$, and as a dominant (resp. recessive) solution it is unique up to a constant multiple, \cite{wtp}, \cite{ph1}. \\

\noindent In this paragraph we fix a pair of dominant/recessive solutions of \eqref{ap1}, say, $A_n$ and $U_n$ respectively. Let $L>0$. Then there is a sequence of reals of the form $V_n/A_n$, where $V_n$ is a solution of \eqref{ap1} such that $V_n/A_n \to L$ as $n \to \infty$. Indeed, choose $V_n$ by setting $V_n = U_n + LA_n$. Hence, for example, there exists a solution $V_n$ of \eqref{ap1} such that $$V_n/A_n \to \zeta(5) \quad n \to \infty,$$ or another solution $W_n$ such that $$W_n/A_n \to \zeta(7) \quad n \to \infty, {\rm etc.}$$ but the terms of $V_n, W_n, {\rm etc.}$ are not necessarily all rational. In addition, for a given real $L>0$ and any $\gamma > 0$, the solution $V_n \equiv B_n + \gamma U_n $ is such that $V_n/A_n \to \zeta(3),$ as $n \to \infty$.\\

\subsection{On the quadratic irrationality of $\zeta(3)$}
\label{qi}

\noindent Another question is whether $\zeta(3)$ is itself algebraic of degree 2 over $\mathbb{Q}$? Although we do not answer this question either way, we present an apparently tractable equivalent formulation which may shed some light on this question. The method is sufficiently general so as to show that given any number known to be irrational by applying an Ap\'{e}ry-type argument on a three term recurrence relation or issuing from a continued fraction expansion, the statement that it is a quadratic irrational is equivalent to a statement about rational valued principal solutions of a corresponding disconjugate four-term recurrence relation.\\

\noindent We proceed first by showing that solutions of a linear three-term recurrence relation can be used to generate a basis for a corresponding four-term linear recurrence relation. The analogous result for differential equations is sufficiently well-known and old, see e.g., Ince~\cite{ince}. Our corresponding result, Theorem~\ref{th2} below, appears to be new in the general case. As a consequence, the quantities $A_n, B_n$ defined in \eqref{an}, \eqref{bn} can be used to generate a basis for a new recurrence relation of order one higher than the original one \eqref{ap1} considered by Ap\'{e}ry.\\

\noindent Given any three-term recurrence relation in general form 
\begin{equation}
\label{g3trr}
p_ny_{n+1} + q_{n}y_{n-1}=r_ny_n,\quad n\geq 1,
\end{equation}
\noindent the mere assumption that $p_nq_n\neq 0$ for all $n$, enables one to transform \eqref{g3trr} into the self-adjoint form \eqref{eqq} below by means of the substitution $c_n=c_{n-1}p_n/q_n$, $c_0$ given, and $b_n=c_nr_n/p_n$. Hence, for simplicity and ease of exposition we assume that the recurrence relation is already in self-adjoint form, and there is no loss of generality in assuming this. We maintain the use of the symbols $A_n, B_n$ for the solutions under consideration for motivational purposes.\\

\begin{theorem}\label{th2} Let $A_n, B_n$ generate a basis for the solution space of the three term recurrence relation \eqref{eqq}
\begin{equation}
\label{eqq}
c_ny_{n+1}+c_{n-1}y_{n-1} - b_ny_n=0, \quad n\geq 1,
\end{equation}
where $c_n \neq 0$, $b_n\neq 0$ for every $n$, and $b_n, c_n \in \mathbb{R}$. Then the quantities $x_{n-1}\equiv A_nA_{n-1}$, $y_{n-1}\equiv B_nB_{n-1}$, $z_{n-1}\equiv A_nB_{n-1}+A_{n-1}B_n$ form a basis for the solution space of the four-term recurrence relation
\begin{eqnarray}
\label{eqqqq}
&& c_{n+2}c_{n+1}^2b_{n}\, z_{n+2} + (b_{n}c_{n+1}^3-b_nb_{n+1}b_{n+2}c_{n+1})\, z_{n+1} + \nonumber\\
 && (b_nb_{n+1}b_{n+2}c_{n} - b_{n+2}c_{n}^3)\,z_n - c_{n-1}c_{n}^2b_{n+2}\, z_{n-1} = 0,\quad n\geq 1,
\end{eqnarray}
\end{theorem}

\begin{proof}(Theorem~\ref{th2}) Direct verification using repeated applications of \eqref{eqq} and simplification, we omit the details. The linear independence can be proved using Wronskians, see below (and see Hartman \cite{ph1} but where in Proposition 2.7 on p. 8 of \cite{ph1}, the reader should replace $a$ by $\alpha$).
\end{proof}
\noindent The Wronskian of the three solutions $x_n = A_{n+1}A_{n}$, $y_n=B_{n+1}B_{n}$, $z_n = A_{n+1}B_{n}+A_{n}B_{n+1}$ of \eqref{eqqqq} arising from the two independent solutions $A_n, B_n$ of the three-term recurrence relation \eqref{eqq} is given by the determinant of the matrix [\cite{ph2}, p.310],

 $$\left[ \begin {array}{ccc} A_{{n+1}}A_{{n}}&B_{{n+1}}B_{{n}}&A_{{n+1}
}B_{{n}}+A_{{n}}B_{{n+1}}\\\noalign{\medskip}A_{{n+2}}A_{{n+1}}&B_{{n+
2}}B_{{n+1}}&A_{{n+2}}B_{{n+1}}+A_{{n+1}}B_{{n+2}}\\\noalign{\medskip}
A_{{n+3}}A_{{n+2}}&B_{{n+3}}B_{{n+2}}&A_{{n+3}}B_{{n+2}}+A_{{n+2}}B_{{
n+3}}\end {array} \right]  $$

\noindent which, after the usual iterations (or see [\cite{ph1}, Prop.2.7]) reduces to the expression:
\begin{equation}
\label{wron}
{\frac {b _{n+2} b _{n+1} c _{n-1} ^{3} \left( A_{{n}}B_{{n-1}}-B_{{n}}A_{{n-1}}
 \right) ^{3}}{c_{n}c_{n+2}c_{n+1}^{3}}}.
\end{equation}
\noindent We apply Theorem~\ref{th2} to the questions at hand, although it is likely there are more numerous applications elsewhere. Thus, the following corollary (stated as a theorem) is immediate.

\begin{theorem}\label{th33} Let $A_n, B_n$ be the Ap\'{e}ry sequences define above in \eqref{an}, \eqref{bn} and consider the corresponding three-term recurrence relation \eqref{ap1} where, for our purposes, $c_n=(n+1)^3, b_n = 34n^3+51n^2+27n+5$. Then the four-term recurrence relation 
\begin{eqnarray}
\label{eqqq}
&& \left( n+3 \right) ^{3} \left( n+2 \right) ^{6} \left( 2\,n+1
 \right)  \left( 17\,{n}^{2}+17\,n+5 \right) z_{{n+2}} \nonumber \\ 
&&- \left( 2\,n+1
 \right)  \left( 17\,{n}^{2}+17\,n+5 \right)  ( 1155\,{n}^{6}+
13860\,{n}^{5}+68535\,{n}^{4} \nonumber \\ 
&&+ 178680\,{n}^{3} 
 + 259059\,{n}^{2}+198156\,
n+62531)  \left( n+2 \right) ^{3}z_{{n+1}} \nonumber \\ 
&&+ \left( 2\,n+5
 \right)  \left( 17\,{n}^{2}+85\,n+107 \right)  ( 1155\,{n}^{6}+
6930\,{n}^{5}+16560\,{n}^{4}\nonumber \\ 
&& +20040\,{n}^{3} + 12954\,{n}^{2}+4308\,n+584)  \left( n+1 \right) ^{3}z_{{n}} \nonumber \\ 
&&- \left( n+1 \right) ^{6}{n}^{
3} \left( 2\,n+5 \right)  \left( 17\,{n}^{2}+85\,n+107 \right) z_{{n-1
}} = 0,
\end{eqnarray}
admits each of the three products $x_{n-1}\equiv A_nA_{n-1}$, $y_{n-1}\equiv B_nB_{n-1}$, and $z_{n-1}\equiv A_nB_{n-1}+A_{n-1}B_n$ as a solution, and the resulting set is a basis for the solution space of \eqref{eqqq}.
\end{theorem}

\noindent The calculation of the Wronskian in the case of \eqref{eqqq} is an now easy matter (see \eqref{wron}). In the case of our three solutions of \eqref{eqqq}, namely $x_{n}, y_{n}, z_{n}$ defined in Theorem~\ref{th2}, the Wronskian is given by 
$${\frac { \left( 2\,n+3 \right)  \left( 2\,n+5 \right)  \left( 17\,{n}^
{2}+51\,n+39 \right)  \left( 17\,{n}^{2}+85\,n+107 \right) {n}^{9}
 \left( A_{{n}}B_{{n-1}}-B_{{n}}A_{{n-1}} \right) ^{3}}{ \left( n+1
 \right) ^{3} \left( n+2 \right) ^{9} \left( n+3 \right) ^{3}}}$$
The non-vanishing of the determinant  for every $n$ is also clear. The counterpart to \eqref{eq6} in this higher order setting is
$$R_n\, \left( A_{{n}}B_{{n-1}}-B_{{n}}A_{{n-1}} \right) ^{3} = L_n \det W(x,y,z)(0),$$
where $W(x,y,z)(0) = -62595/64$,
$$R_n \equiv {\frac { \left( 2\,n+3 \right)  \left( 2\,n+5 \right)  \left( 17\,{n}^
{2}+51\,n+39 \right)  \left( 17\,{n}^{2}+85\,n+107 \right) {n}^{9}
 }{ \left( n+1
 \right) ^{3} \left( n+2 \right) ^{9} \left( n+3 \right) ^{3}}},$$
and 
$$L_n \equiv \prod_{m=1}^{n}{\frac{m^3(m+1)^6(2m+5)(17m^2+85m+107)}{(m+2)^6(m+3)^3(2m+1)(17m^2+17m+5)}}.$$

\noindent Recall that $A_nA_{n-1}$ is a solution of \eqref{eqqq}, that $A_0=1$, $A_1=5$, and $A_n > 0$ for every $n > 1$.

\begin{theorem}\label{th3} $\zeta(3)$ is algebraic of degree two over $\mathbb{Q}$ if and only if \eqref{eqqq} has a non-trivial rational valued solution $S_n$ (i.e., $S_n$ is rational for every $n\geq 1$), with 
\begin{equation}
\label{sn0}\frac{S_n}{A_nA_{n-1}} \to 0, \quad n \to \infty.
\end{equation}
\end{theorem}

\begin{proof}(Theorem~\ref{th3}) (Sufficiency) Since $A_nA_{n-1}$, $B_nB_{n-1}$ and $A_{n-1}B_n+A_nB_{n-1}$ are linearly independent we have 
\begin{equation}\label{sn}S_n = aA_nA_{n-1}+bB_nB_{n-1}+c(A_{n-1}B_n+A_nB_{n-1}),
\end{equation} 
for some $a,b,c \in \mathbb{R}$, not all zero. Since $S_n$ is rational valued for all $n$ by hypothesis, the repeated substitutions $n=1,2,3$ in the above display yield a system of three equations in the unknowns $a,b,c.$ Since all the values involved are rational numbers, the same is true of this unique set of $a,b,c$. \\

With this set of $a,b,c$ we note that, for every $n \geq 1$, \begin{equation}\label{sn1}\frac{S_n}{A_nA_{n-1}} = a\,+ b\, \frac{B_nB_{n-1}}{A_nA_{n-1}}+c\, \left (\frac{B_n}{A_n}+\frac{B_{n-1}}{A_{n-1}}\right ).\end{equation} But from Ap\'{e}ry's work \cite{rap} (or \cite{rm}, \cite{vdp}) we know that $B_n/A_n \to \zeta(3)$ as $n \to \infty$. Using this information in passing to the limit we have that \begin{equation}\label{sn2}\frac{S_n}{A_nA_{n-1}} \to a \, + b\,\zeta(3)^2 + 2c\, \zeta(3), \quad n \to \infty.\end{equation} 
The possibility that $b=0$ is excluded by the fact that $\zeta(3)$ is irrational. Thus, $b \neq 0$ and so $\zeta(3)$ is algebraic of degree 2 over $\mathbb{Q}$. \\

\noindent Conversely, assume that $\zeta(3)$ is algebraic of degree 2 over $\mathbb{Q}$. Then, there exists rational constants $a,b,c$ with $b\neq 0$ such that $b\,\zeta(3)^2 + 2c\, \zeta(3) + \,a =0.$ For this choice of $a,b,c,$ consider the solution of \eqref{eqqq} defined by \eqref{sn}. Since this $S_n$ clearly satisfies \eqref{sn1}, and \eqref{sn2} by construction, the limiting result \eqref{sn0} follows.
\end{proof}

\noindent The proof of Theorem~\ref{th3} is capable of much greater generality. Combined with Theorem~\ref{th2} and minor changes in the argument of the previous theorem one can easily prove

\begin{theorem}\label{th4} Let $c_n, b_n \in \mathbb{R}$, $c_n \neq 0$, $b_n\neq 0$, for every $n$. Let $A_n, B_n$ be two independent rational valued solutions of \eqref{eqq} such that $$\lim_{n\to \infty}\frac{B_n}{A_n} = L,$$
where $L$ is irrational. Then $L$ is algebraic of degree two over $\mathbb{Q}$ if and only if \eqref{eqqqq} has a non-trivial rational valued solution $S_n$ such that \eqref{sn0} holds.
\end{theorem}

\begin{example} \label{exxx}{\rm Consider the Fibonacci sequence $F_n$ defined by the self-adjoint three term recurrence relation \eqref{eqq} with $b_n=c_n=(-1)^n$ for all $n\geq 1$, and $c_0=1.$ Then it is easy to see that the solutions defined by $A_0=0$, $A_1=1$, $B_0=1$ and $B_1=1$ are given by $A_n \equiv F_{n}$, $B_n = F_{n+1}$ are two linearly independent solutions of the Fibonacci relation such that $B_n/A_n \to L$, where $L = (1+\sqrt{5})/2$ is already known to be irrational (it is not necessary that $L$ be irrational as a result of the actual approach to the limit).\\ 

\noindent According to Theorem~\ref{th2} the quantities $F_nF_{n-1}$, $F_{n+1}F_n$ and $F_{n}^2 + F_{n-1}F_{n+1}$ satisfy the four-term recurrence relation $$z_{n+2}-2z_{n+1}-2z_{n}+z_{n-1}=0.$$ Note that the solution $S_n$ defined by 
$$S_n = - F_nF_{n-1} + F_{n+1}F_n - (F_{n}^2 + F_{n-1}F_{n+1})/2,$$ is a nontrivial rational valued solution of this four term recurrence relation such that $S_n/F_nF_{n-1} \to L^2-L-1 = 0$ as $n \to \infty$. Hence, by Theorem~\ref{th2}, $L$ is algebraic of degree 2. }
\end{example}

\vskip0.25in

\subsection{Discussion} In the language of the theory of disconjugate difference equations a {\it special} solution like $S_n$ in Theorem~\ref{th3}, if it exists, is a 2nd principal solution of \eqref{eqqq}. In the case of a disconjugate four-term linear recurrence relation with a positive leading term (such as ours, \eqref{eqqq}), a ${\rm k^{th}}$ principal solution $u_{k, n}$ is characterized by the existence of limits of the form $$\frac{u_{k-1, n}}{u_{k, n}} \to 0, \quad n \to \infty,$$ for $k=1,2,3$. A first principal solution $u_{0,n}$ is unique up to a multiplicative constant, when it exists. For example, in the case of $\zeta(3)$, \eqref{eq8} gives us that quotients of solutions of disconjugate linear recurrence relations always have limits at infinity (and they are allowed to be infinite). In the case of \eqref{eqqq} this is easy to see since we know the basis explicitly. For example, the limit of the two solutions $u_n = B_nB_{n-1}$ and $v_n=2\zeta(3)A_nA_{n-1}-(A_{n}B_{n-1}+A_{n-1}B_n)$ of \eqref{eqqq} exists at infinity, and $\lim_{n\to \infty}u_n/v_n = + \infty$.  On the other hand, the solution $v_n$ just defined and $w_n = A_nA_{n-1}+A_nB_{n-1}+A_{n-1}B_n$ are such that $\lim_{n\to \infty}v_n/w_n = 0.$ In the case of disconjugate or non-oscillatory difference equations (or recurrence relations) such principal solutions always exist, see Hartman [\cite{ph1}, Section 8], [\cite{ph2}, Appendix A] for basic discussions on these and related matters.\\

\noindent

\section{A criterion for algebraic irrationality and transcendence}

\noindent In this section we show how the results of the previous sections may be regarded as special cases of a more general application. Indeed, assume that we have an irrational number $L$ whose rational approximations are derived either by means of an Ap\'{e}ry type argument on a three term recurrence relation, or perhaps via a continued fraction expansion of $L$. Basically, we show that if $L$ is not algebraic of degree less than or equal to $(m-1)$, then $L$ is algebraic of degree $m$ over $\mathbb{Q}$ if and only if there exists a disconjugate $(m+2)-$term linear recurrence relation having a non-trivial rational valued principal solution of a specific type. A special case of the result to follow is to be found in Theorem~\ref{th4} above. This of course, also leads to a necessary and sufficient condition for the transcendence of such numbers. We outline herewith the construction of this special recurrence relation pointing out first two important special cases as motivation: The first case is to be found in Theorem~\ref{th4} as alluded to above. The second case is a ``degree 3" version of Theorem~\ref{th4} which we describe next.\\

\noindent Associated to \eqref{eqq} is a higher order analog of Theorem~\ref{th2}. Consider the 5-term recurrence relation
\begin{equation}\label{eqqq5}p_n\,z_{n+3}-q_n\,z_{n+2}-r_n\,z_{n+1}-s_n\,z_{n}-t_n\,z_{n-1}=0, \quad n\geq 1,
\end{equation}
where the leading term $p_n\neq 0$ for all $n$, and\\

\noindent $p_n = c_{{n+4}}{c_{{n+3}}}^{2}{c_{{n+2}}}^{3}b_{{n+1}} \left( {c_{{n}}}^{2} -b_{{n+1}}b_{{
n}}\right),$

\noindent $q_n = -b_{{n+2}}b_{{n+1}}{c_{{n+3}}}^{3}{c_{{n+2}}}^{2}{c_{{n}}}^{2}+{b_{{n+
1}}}^{2}c_{{n+3}}{c_{{n+2}}}^{4}b_{{n+4}}b_{{n}}- \\ 
b_{{n+2}}{b_{{n+1}}}^{2}b_{{n+3}}c_{{n+3}}{c_{{n+2}}}^{2}b_{{n+4}}b_{{n}}+b_{{n+2}}b_{{n+1}
}b_{{n+3}}c_{{n+3}}{c_{{n+2}}}^{2}b_{{n+4}}{c_{{n}}}^{2}- \\ 
b_{{n+1}}c_{{n+3}}{c_{{n+2}}}^{4}b_{{n+4}}{c_{{n}}}^{2}+b_{{n+2}}{b_{{n+1}}}^{2}{c_
{{n+3}}}^{3}{c_{{n+2}}}^{2}b_{{n}}$\\

\noindent $r_n = -b_{{n+2}}b_{{n+1}}{b_{{n+3}}}^{2}c_{{n+2}}c_{{n+1}}b_{{n+4}}{c_{{n}}}
^{2}-{b_{{n+1}}}^{2}{c_{{n+2}}}^{3}c_{{n+1}}b_{{n+4}}b_{{n}}b_{{n+3}}+ \\
b_{{n+1}}b_{{n+3}}c_{{n+2}}{c_{{n+1}}}^{3}{c_{{n+3}}}^{2}b_{{n}}+b_{{n+2}}b_{{n+1}}b_{{n+3}}c_{{n+2}}c_{{n+1}}{c_{{n+3}}}^{2}{c_{{n}}}^{2}- \\ 
b_{{n+1}}c_{{n+2}}{c_{{n+1}}}^{3}b_{{n+4}}b_{{n}}{b_{{n+3}}}^{2}+{b_{{n
+1}}}^{2}{c_{{n+2}}}^{3}c_{{n+1}}{c_{{n+3}}}^{2}b_{{n}} + \\ 
b_{{n+1}}b_{{n+3}}{c_{{n+2}}}^{3}c_{{n+1}}b_{{n+4}}{c_{{n}}}^{2}+{b_{{n+3}}}^{2}c_{{
n+2}}{c_{{n+1}}}^{3}b_{{n+4}}{c_{{n}}}^{2}- \\ 
b_{{n+1}}{c_{{n+2}}}^{3}c_{{n+1}}{c_{{n+3}}}^{2}{c_{{n}}}^{2}-b_{{n+2}}{b_{{n+1}}}^{2}b_{{n+3}}c_
{{n+2}}c_{{n+1}}{c_{{n+3}}}^{2}b_{{n}}+ \\ 
b_{{n+2}}{b_{{n+1}}}^{2}c_{{n+2}}c_{{n+1}}b_{{n+4}}b_{{n}}{b_{{n+3}}}^{2}-b_{{n+3}}c_{{n+2}}{c_{{n+1}
}}^{3}{c_{{n+3}}}^{2}{c_{{n}}}^{2},$\\

\noindent $s_n=- b_{{n+3}}{c_{{n+1}}}^{4}c_{{n}}{c_{{n+3}}}^{2}b_{{n}}+b_{{n+2}}{c_{{n
+1}}}^{2}{c_{{n}}}^{3}b_{{n+4}}{b_{{n+3}}}^{2}+ \\ 
{b_{{n+3}}}^{2}{c_{{n+1}}}^{4}c_{{n}}b_{{n+4}}b_{{n}}+ b_{{n+2}}b_{{n+1}}b_{{n+3}}{c_{{n+1}}}^
{2}c_{{n}}{c_{{n+3}}}^{2}b_{{n}}- \\ 
b_{{n+2}}b_{{n+1}}{b_{{n+3}}}^{2}{c_{{n+1}}}^{2}c_{{n}}b_{{n+4}}b_{{n}}- b_{{n+2}}{c_{{n+1}}}^{2}{c_{{n}}}^{
3}{c_{{n+3}}}^{2}b_{{n+3}},$\\

\noindent $t_n={b_{{n+3}}}^{2}b_{{n+4}}{c_{{n+1}}}^{3}{c_{{n}}}^{2}c_{{n-1}}-b_{{n+3}
}{c_{{n+3}}}^{2}{c_{{n+1}}}^{3}{c_{{n}}}^{2}c_{{n-1}}.$\\

\noindent Note that the hypothesis, $p_n\neq 0$ for all $n$, is equivalent to $t_n\neq 0$ for all $n$. Then for any given pair of linearly independent solutions $A_n, B_n$ of \eqref{eqq} the sequences $A_{n+1}A_nA_{n-1}$, $B_{n+1}B_nB_{n-1}$, $A_{n+1}A_nB_{n-1}+A_{n+1}B_nA_{n-1}+B_{n+1}A_nA_{n-1}$, and $B_{n+1}B_nA_{n-1}+B_{n+1}A_nB_{n-1}+A_{n+1}B_nB_{n-1}$, form a linearly independent set of solutions for \eqref{eqqq5}. Given that we know how to test for degree 2 irrationality of limits $L$ via Theorem~\ref{th4}, we can formulate an analogous result for degree 3 irrationality next.\\

\noindent {\bf Note:} In the sequel, we always assume that the $A_n, B_n$ in question are positive for all $n$ (as they arise from a disconjugate equation \eqref{eqq}). There is no loss of generality in assuming this since the proofs involve limiting arguments. Also, unless otherwise specified we assume that $L\neq0$.\\

\begin{theorem}\label{th5} Let $c_n, b_n \in \mathbb{R}$ and $p_n \neq 0$ for all $n$. Let $A_n, B_n$ be two independent rational valued solutions of \eqref{eqq} such that $$\lim_{n\to \infty}\frac{B_n}{A_n} = L,$$
where $L$ is irrational and $L$ is not algebraic of degree 2. Then $L$ is algebraic of degree three over $\mathbb{Q}$ if and only if \eqref{eqqq5} has a non-trivial rational valued solution $S_n$ such that 
\begin{equation}\label{sn005} \frac{S_n}{A_{n+1}A_{n}A_{n-1}} \to 0, \quad n\to \infty.
\end{equation}
\end{theorem}
\noindent The proof is similar to that of Theorem~\ref{th3} and so is omitted. \\

\begin{remark}{\rm A re-examination of the proof of Theorem~\ref{th3} which serves as a template for all other such proofs to follow shows that the tacit assumptions on $L$ can be waived to some extent. The previous result may then be re-formulated as follows.\\
}\end{remark}

\noindent Let $S_n$, a solution of \eqref{eqqq5}, have the basis representation 
\begin{eqnarray}
\label{basis}
&& S_n =a_3\,A_{n+1}A_nA_{n-1} +a_0\, B_{n+1}B_nB_{n-1}+ \nonumber \\ 
&& a_2\,(A_{n+1}A_nB_{n-1}+A_{n+1}B_nA_{n-1}+B_{n+1}A_nA_{n-1}) + \nonumber \\ 
&& a_1\,(B_{n+1}B_nA_{n-1}+B_{n+1}A_nB_{n-1}+A_{n+1}B_nB_{n-1}),
\end{eqnarray}
where $a_i \in \mathbb{R}$ and the subscript $i$ in $a_i$ for the basis coordinates is determined by counting the number of $A$'s in the basis vector immediately following it. \\

\begin{theorem}\label{deg3} Let $c_n, b_n \in \mathbb{R}$ and $p_n \neq 0$ for all $n$. Let $A_n, B_n$ be two independent rational valued solutions of \eqref{eqq} such that $$\lim_{n\to \infty}\frac{B_n}{A_n} = L.$$
Then $L$ is algebraic of degree at most 3 if and only if there exists a non-trivial rational valued solution $S_n$ of \eqref{eqqq5}  satisfying \eqref{sn005}.
\end{theorem}

\begin{proof}(Theorem~\ref{deg3}) Idea: Using \eqref{basis} we see that since $S_n$ is rational, then so are the $a_i$, $0\leq i \leq 3$, not all of which are zero. Next, as $n \to \infty$,
$$\frac{S_n}{A_{n+1}A_nA_{n-1} }  \to a_0L^3+ 3a_1L^2+3a_2L+a_3,$$ and so $L$ is algebraic of degree no greater than 3. Conversely, let $L$ be algebraic of degree no greater than 3 and let $p(x) = a_0x^3+ 3a_1x^2+3a_2x+a_3$ be its defining polynomial where not all $a_i$ are zero. Then choosing the solution $S_n$ of \eqref{eqqq5} in the form \eqref{basis} with the same quantities $a_i$ that appear as the coefficients of $p$, we see that since $p(L)=0$, \eqref{sn005} is satisfied.
\end{proof}

\begin{remark}{\rm In order to improve on Theorem~\ref{th5} we need to add more to the solution $S_n$ appearing therein. For example, it is easy to see that under the same basic conditions on the $A_n, B_n$,  if there exists a non-trivial rational valued solution $S_n$ of \eqref{eqqq5} with $a_0 \neq 0$ satisfying \eqref{sn005}, then $L$ is algebraic of degree no greater than 3. On the other hand, if $L$ is algebraic of degree 3, then there exists a non-trivial rational valued solution $S_n$ of \eqref{eqqq5} with $a_0 \neq 0$ satisfying \eqref{sn005}.
}\end{remark}

\noindent Theorem~\ref{th4} and Theorem~\ref{th5} give us an idea on how to proceed next. In essence, we now have some way of determining whether or not the limit $L$ is algebraic of degree 3 based on the fact that it is not algebraic of lower degree. The general result is similar, but first we describe the construction of the required linear higher order recurrence relations. In order to test whether the limit $L$ in Theorem~\ref{th5} is algebraic of degree $m$, $m \geq 2$, we will require a linear recurrence relation containing $(m+2)-$terms or equivalently an $(m+1)-$th order linear difference equation (an equation involving ``finite differences" in the traditional sense). This new equation is found from a prior knowledge of the kernel of the associated operator. \\

\noindent As usual we let $A_n, B_n$ be two linearly independent solutions of \eqref{eqq}. We seek a homogeneous linear $(m+2)-$term recurrence relation whose basis (consisting of $(m+1)$ terms) is described as follows: Two basic elements are given by 
$$A_{n+m-2}A_{n+m-1}\ldots A_nA_{n-1}$$
along with a corresponding term with all these $A$'s replaced by $B$'s. To each given $k$, $0<k<m$, we associate a sum of products of terms of the form 
$$\sum A_{n+m-2}A_{n+m-1}\ldots B_{n+m-i}\ldots B_{n+m-j}  \ldots A_nA_{n-1}$$
where this sum contains exactly $\binom{m}{k}$ distinct terms. Each summand is obtained by enumerating all the possible ways of choosing $k-$terms out of the full product of $A$'s and replacing each such $A$ by a $B$ (while keeping the subscripts intact).\\

\noindent For example, if $m=4$ and $k=2$ there is a such a sum of $6=\binom{4}{2}$ terms, the totality of which is of the form \\ 
\begin{eqnarray*}
&&x_{3}^{(n-1)} = A_{n+2}A_{n+1}B_{n}B_{n-1}+A_{n+2}B_{n+1}B_{n}A_{n-1}+B_{n+2}B_{n+1}A_{n}A_{n-1}+ \\ 
&&B_{n+2}A_{n+1}B_{n}A_{n-1}+A_{n+2}B_{n+1}A_{n}B_{n-1}+B_{n+2}A_{n+1}A_{n}B_{n-1}.\\
\end{eqnarray*}

\noindent The collection of all such ``sums of products" as $k$ varies from 0 to m gives us a collection of $(m+1)$ elements denoted by $x_{1}^{(n-1)}, x_{2}^{(n-1)},\ldots,x_{m+1}^{(n-1)}$. That this specific set of elements is a linearly independent set may depend on the nature of the interaction of the $a_n, b_n$ in \eqref{eqq} as we saw above (e.g., $p_n\neq 0$ in \eqref{eqqq5}). At any rate, since every solution $z_{n-1}$ of this new recurrence relation must be a linear combination of our $x_{i}^{(n-1)}$, it is easy to see that the compatibility relation is obtained by setting the determinant of the matrix 
$$\left[ \begin {array}{ccccc}
z_{n-1}&x_{1}^{(n-1)}&x_{2}^{(n-1)}&\ldots& x_{m+1}^{(n-1)}\\\noalign{\medskip}
z_{n}&x_{1}^{(n)}&x_{2}^{(n)}&\ldots& x_{m+1}^{(n)}\\\noalign{\medskip}
\ldots&\ldots&\ldots&\ldots& \ldots\\\noalign{\medskip}
\ldots&\ldots&\ldots&\ldots& \ldots\\\noalign{\medskip}
z_{n+m}&x_{1}^{(n+m)}&x_{2}^{(n+m)}&\ldots& x_{m+1}^{(n+m)}\\\noalign{\medskip}
\end {array} \right],$$
equal to zero, for every $n$. This and the repeated use of the recurrence relation \eqref{eqq} gives the required $(m+2)-$term recurrence relation of which \eqref{eqqq5} and \eqref{eqqqq} are but special cases.\\ \\

\noindent  {\bf Note:} In the sequel we always assume that the set consisting of the ``sums of products" described above is a linearly independent set of $(m+1)$ elements. This is equivalent to various conditions to be imposed upon the coefficient of the leading and trailing terms of the ensuing $(m+2)-$term recurrence relation whose construction is presented above.\\ \\

\begin{theorem}\label{th6} Let $c_n, b_n \in \mathbb{R}$, $c_n\neq 0$ and $b_n\neq 0$ in addition to other conditions enunciated in the note above. Let $m\geq 3$. Consider two independent rational valued solutions $A_n, B_n$ of \eqref{eqq} such that $$\lim_{n\to \infty}\frac{B_n}{A_n} = L,$$
where $L$ is not algebraic of degree less than or equal to $m-1$. Then $L$ is algebraic of degree $m$ over $\mathbb{Q}$, if and only if the $(m+2)-$term linear recurrence relation described above has a non-trivial rational valued solution $S_n$ such that 
\begin{equation}\label{sn0005} \frac{S_n}{A_{n+m-2}\ldots A_{n+1}A_{n}A_{n-1}} \to 0, \quad n\to \infty.
\end{equation}
\end{theorem}
\vskip0.25in
\noindent An analog of Theorem~\ref{deg3} can also be formulated, perhaps easier to use in practice.\\

\begin{theorem}\label{degn} Let $c_n, b_n \in \mathbb{R}$, $c_n\neq 0$ and $b_n\neq 0$ in addition to other conditions enunciated in the note above. Let $m\geq 3$. Let $A_n, B_n$ be two independent rational valued solutions of \eqref{eqq} such that $$\lim_{n\to \infty}\frac{B_n}{A_n} = L.$$
Then $L$ is algebraic of degree at most $m$ over $\mathbb{Q}$, if and only if the $(m+2)-$term linear recurrence relation described above has a non-trivial rational valued solution $S_n$ such that 
\begin{equation}\label{sn0005} \frac{S_n}{A_{n+m-2}\ldots A_{n+1}A_{n}A_{n-1}} \to 0, \quad n\to \infty.
\end{equation}
\end{theorem}

\begin{remark}{\rm  Since the condition in Theorem~\ref{degn} puts a bound on the degree $m$ of algebraic irrationality over $\mathbb{Q}$ it also gives an equivalent criterion for the transcendence of numbers $L$ whose limits are found by using quotients of solutions of three-term recurrence relations. In particular, associated to the special number $\zeta(3)$ is an infinite sequence of specific linear recurrence relations of every order, as constructed above, involving sums of products of both sets of Ap\'{e}ry numbers $A_n, B_n$. The transcendence of $\zeta(3)$ is then equivalent to the statement that none of the infinite number of (disconjugate) recurrence relations constructed has a nontrivial rational valued principal solution of the type described.}
\end{remark}

\begin{example}{\rm In this final example we interpret Ap\'{e}ry's construction \cite{rap}, for the irrationality of $\zeta(2)$ in the context of the non-existence of rational valued solutions of recurrence relations with predetermined asymptotics. Recall that Ap\'{e}ry's three term recurrence relation for the proof of the irrationality of $\zeta(2)$ is given by \cite{rap}
$$(n+1)^2y_{n+1} - n^2y_{n-1} = (11n^2+11n+3)y_{n},\quad n\geq 1.$$
In order to apply Theorem~\ref{th2} we need to express this equation in self-adjoint form; that is we simply multiply both sides by $(-1)^n$ resulting in the equivalent equation \eqref{eqq} with $c_n = (-1)^n(n+1)^2$, $b_n=(-1)^n(11n^2+11n+3)$. The Ap\'{e}ry solutions of this equation (e.g., \cite{vdp}) are given by
\begin{equation}\label{an2}
A_{n}^{\prime} = \sum_{k=0}^{n}{\binom{n}{k}}^2{\binom{n+k}{k}},
\end{equation}
and
\begin{equation}\label{bn2}
B_{n}^{\prime} = \sum_{k=0}^{n}{\binom{n}{k}}^2{\binom{n+k}{k}}\left\{2\sum_{m=1}^{n} \frac{(-1)^{m-1}}{m^2} + \sum_{m=1}^{k} \frac{(-1)^{n+m-1}}{m^2\binom{n}{m}\binom{n+m}{m}} \right \}.
\end{equation}
From Ap\'{e}ry's work it is known that these solutions have the property that $$B_{n}^{\prime}/A_{n}^{\prime} \to \zeta(2)$$ as $n \to \infty$, and we already know that $\zeta(2)$ is irrational. It follows from the above considerations that the four term recurrence relation
\begin{eqnarray*} p_nz_{n+2}+q_nz_{n+1}+r_nz_n+s_nz_{n-1}=0,
\end{eqnarray*}
where \\ \\
$p_n= \left( n+3 \right) ^{2} \left( n+2 \right) ^{4} \left( 11\,{n}^{2}+11
\,n+3 \right), $\\
$q_n = - \left( 11\,{n}^{2}+11\,n+3 \right)  \left( 122\,{n}^{4}+976\,{n}^{3}
+2873\,{n}^{2}+3684\,n+1741 \right)  \left( n+2 \right) ^{2},$\\
$r_n = - \left( 11\,{n}^{2}+55\,n+69 \right)  \left( 122\,{n}^{4}+488\,{n}^{3
}+677\,{n}^{2}+378\,n+76 \right)  \left( n+1 \right) ^{2},$\\
$s_n = \left( n+1 \right) ^{4}{n}^{2} \left( 11\,{n}^{2}+55\,n+69 \right),$\\ \\
and whose basis is given by the three elements $A_{n}^{\prime}A_{n-1}^{\prime}$, $B_{n}^{\prime}B_{n-1}^{\prime}$ and $A_{n}^{\prime}B_{n-1}^{\prime}+B_{n}^{\prime}A_{n-1}^{\prime}$ cannot have a non-trivial rational valued solution $S_n$ satisfying $$\frac{S_n}{A_{n}^{\prime}A_{n-1}^{\prime}} \to 0, \quad n\to \infty.$$ But we also know that $\zeta(2)$ is actually transcendental (as it is a rational multiple of $\pi^2$), and so cannot be algebraic of any finite degree. Hence, for each $m$, none of the $(m+2)-$term recurrence relations that can be constructed as described above has a nontrivial rational valued solution satisfying \eqref{sn005}.}
\end{example}

\section{Appendix}

\noindent The following series for the integer powers of $\zeta(3)$ were motivated by the results of the last section. In what follows $A_n, B_n$ are the standard Ap\'{e}ry sequences defined by \eqref{an}, \eqref{bn} above and $b_n$ is the Ap\'{e}ry cubic defined in \eqref{ap1}, that is $b_n = 34n^3+51n^2+27n+5$. Recall that the first series on the following list is \eqref{eq15}, above:
{ \mathversion{bold}
\begin{equation*}
\zeta (3) =  \,6\,\sum_{n=1}^{\infty}\frac{1}{n^3A_nA_{n-1}},
\end{equation*}}

{ \mathversion{bold}
$$\zeta(3)^2 = 6\,\sum _{n=1}^{\infty}\,\frac{{p_{n,1}\,B_n}} {n^{3} (n+1) ^{3}A_{n-1}A_{n}A_{n+1}},$$}
\noindent where $p_{n,1}=b_n$,\\

{ \mathversion{bold}
$$\zeta(3)^3 = 6\,\sum _{n=1}^{\infty}\,{\frac { p_{n,1} \, p_{n,2}\, B_{n}\, }
{{n}^{3}\left( n+1 \right) ^{6}\left( n+2 \right) ^{3}A_{{n-1}}A_{{n}}A_{{n+1}} A_{{n+2}} }},$$}

\noindent where\\
\noindent $p_{n,1} = 1155\,{n}^{6}+6930\,{n}^{5}+16560\,{n}^{4}+20040\,{n}^{3}+12954\,{n}^{2}+4308\,n+584,$ and \\
\noindent $p_{n,2} = B_{{n}}b_n -B_{{n-1}}{n}^{3} $,\\

{ \mathversion{bold}
$$\zeta(3)^4 = 12\,\sum _{n=1}^{\infty}{\frac {p_{n,1}\, p_{n,2}\, p_{n,3}\, B_{{n}}}
{{n}^{3}\left( n+1 \right) ^{9}\left( n+2 \right) ^{6} \left( n+3 \right) ^{3}A_{{n-1}}A_{{n}}A_{{n+1}}A_{{n+2}} A_{{n+3}} 
}},$$}
\vskip0.25in
\noindent{w}here\\
\noindent $p_{n,1}=\left( 2\,n+3 \right)(9809\,{n}^{8}+117708\,{n}^{7}+589827\,{n}^{6}+1600641\,{n}^{5}+2554545\,{n}^{4}+$ \\ $+ 2441061\,{n}^{3}+1362947\,{n}^{2}+411198\,n+52020)$\\

\noindent $p_{n,2} = B_{{n}} b_n -B_{{n-1}}{n}^{3}$,\\

\noindent $p_{n,3}= p_{n,3,1}\, B_{{n}} - p_{n,3,2}\,B_{{n-1}}\, n^3,$ \\
\noindent $p_{n,3,1}= 1155\,{n}^{6}+6930\,{n}^{5}+16560\,{n}^{4}+20040\,{n}^{3}+12954\,{n}^{2}+4308\,n+584,$\\
\noindent $p_{n,3,2}=34\,{n}^{3}+ 153{n}^{2}+231\,{n}+117\,,$\\
 \\

{ \mathversion{bold}
\begin{equation*} \zeta(3)^5=6\,\sum _{n=1}^{\infty}\,\frac {p_{n,1}\,p_{n,2}\,p_{n,3}\,p_{n,4}\,B_{n} }{{n}^{3}\left( n+1 \right) ^{12}\left( n+2 \right) ^{9}\left( n+3 \right) ^{6}\left( n+4 \right) ^{3}\prod_{i=n-1}^{n+4}A_i\,},
\end{equation*} } \\
where\\

\noindent $p_{n,1}=1332869\,{n}^{12}+  31988856\,{n}^{11}+342113817\,{n}^{10}+2150577460\,{n}^{9}+8825260041\,{n}^{8}+ 
24829342992\,{n}^{7}+ 48939099945\,{n}^{6}+ 67836980844\,{n}^{5}+ 65389823136\,{n}^{4}+42618151360\,{n}^{3}+ 17812032480\,{n}^{2}+4300387200\,n+456205824$,\\

\noindent$p_{n,2}=b_{n}B_{{n}} -B_{{n-1}}{n}^{3},$\\

\noindent $p_{n,3}= p_{n,3,1}\, B_{{n}} - p_{n,3,2}\, B_{{n-1}}\,n^3,$\\
\noindent $p_{n,3,1} = 39236\,{n}^{9}+529686\,{n}^{8}+3065556\,{n}^{7}+9941526\,{n}^{6}+19822026\,{
n}^{5}+25091514\,{n}^{4}+  20098154\,{n}^{3}+9822474\,{n}^{2}+2675268\,n+312120,$\\
\noindent and $p_{n,3,2}=1155\,{n}^{6}+ 13860\,{n}^{5} +68535\,{n}^{4}+178680
\,{n}^{3}+ 259059\,{n}^{2}+198156\,{n}+62531\,$ \\

\noindent $p_{n,4}= p_{n,4,1}\,B_{{n}} - p_{n,4,2}B_{{n-1}}\,n^3$\\
\noindent $p_{n,4,1}=1155\,{n}^{6}+6930\,{n}^{5}+16560\,{n}^{4}+20040\,{n}^{3}+12954\,{n}^{2}+4308\,n+584$,\\
\noindent $p_{n,4,2}= 34\,{n}^{3}+153\,{n}^{2}+231\,{n} +117,$ and \\
\noindent etc., \\ 

\noindent{w}ith the series of all higher powers of $\zeta(3)$ being exactly computable, the general term being of the form

{ \mathversion{bold}
\begin{equation*} 
\zeta(3)^m=c\,\sum _{n=1}^{\infty}\,{\frac {B_{n} \prod_{i=1}^{m-1} p_{n,i}} {\prod_{i=0}^{m-1}(n+i)^{3(m-i)}\,\prod_{i=n-1}^{n+m-1}A_i}}\quad\quad m\geq 2,
\end{equation*} } \\
\noindent where $c>0$ is a constant, $p_{n,1}$ is a polynomial in $n$ of degree $3(m-1)$, $p_{n,2}=b_{n}B_{{n}} - B_{{n-1}}{n}^{3},$ and generally, for $i\geq 2$, $p_{n,i}=p_{n,i,1}B_n - p_{n,i,2}B_{n-1}n^3,$ where $p_{n,i,j}$ is a polynomial in $n$ of degree $3(i-j)$, for $j=1,2$, and all polynomials above have integer coefficients.\\ \\

\noindent Akin to these series for $\zeta(3)$ are completely analogous corresponding series for powers of $\zeta(2)$, series such as

$$\zeta(2)^2=5\,\sum _{n=1}^{\infty}\left( -1 \right) ^{n-1}{\frac { \left( 11\,{n}^{2}+11\,n+3 \right) B_{n}^{\prime}
 }{{n}^{2} \left( n+1 \right) ^{2}A_{n-1}^{\prime}A_{n+1}^{\prime}A_{n}^{\prime}}} = \pi^4/36,$$

or,

$$\zeta(2)^3 = 5\,\sum _{n=1}^{\infty} \left( -1 \right) ^{n-1}\,{\frac {p_{n,1}p_{n,2} B_{{n}}  }{{n}^{2} \left( n+1 \right) ^{4}\left( n+2 \right) ^{2}A_{n+1}^{\prime}A_{n+2}^{\prime} A_{n}^{\prime}A_{n-1}^{\prime}}} = \pi^6/216,$$

\noindent where\\ \\
\noindent  $p_{n,1}=\left( 122\,{n}^{4}+488\,{n}^{3}+677\,{n}^
{2}+378\,n+76 \right),$\\ \\
\noindent $p_{n,2}= B_{{n-1}}{n}^{2}+ \left( 11\,{n}^
{2}+11\,n+3 \right) B_{{n}}, $\\ \\ \noindent etc.

\end{document}